\documentclass[12pt]{article}
\usepackage{amssymb,amsfonts,amsthm}
\usepackage{amstext}

\usepackage{amsmath}
\usepackage{MnSymbol}
\usepackage{color}
\usepackage[latin1]{inputenc}
\usepackage[active]{srcltx}
\def\openC{{\rm C\kern-.18cm\vrule width.8pt height 7pt depth-.2pt \kern.18cm}}
\def\openN{{{\rm I}\kern-.16em {\rm N}}}
\def\openR{{{\rm I}\kern-.16em {\rm R}}}
\def\openT{{{\rm T}\kern-.42em {\rm T}}}
\def\openZ{{{\rm Z}\kern-.28em{\rm Z}}}

\def\eop{\hfill\rule{2.5mm}{2.5mm}}
\def\pf{\par\smallbreak\noindent {\bf Proof.} \ }

\newtheorem{thm}{Theorem}[section]

\newtheorem{lem}[thm]{Lemma}
\newtheorem{prop}[thm]{Proposition}

\theoremstyle{definition}

\def\eop{\hfill\rule{2.5mm}{2.5mm}}
\begin{document}

\title{\textbf{Strictly positive definite kernels on a product of spheres} \vspace{-4pt}
\author{\sc
J. C. Guella\thanks{Both authors partially supported by FAPESP, under grants $\#$ 2012/22161-3 and $\#$2014/00277-5 and respectively.}\ {\small and}  V. A. Menegatto }}
\date{}
\maketitle \vspace{-30pt}
\bigskip
\begin{center}
\parbox{13 cm}{{\small For the real, continuous, isotropic and positive definite kernels on a product of spheres, one may consider not only its usual strict positive definiteness but also strict positive definiteness restrict to the points of the product that have distinct components.\ In this paper, we provide a characterization for strict positive definiteness in these two cases, settling all the cases but those in which one of the spheres is a circle.}}
\end{center}

\vspace*{3mm}


\noindent{\bf 2010 Math. Subj. Class.}: 33C50; 33C55; 42A82; 43A35

\thispagestyle{empty}

%
%

\section{Introduction}\label{s1}

Let $S^m$ denote the $m$-dimensional unit sphere in $\mathbb{R}^{m+1}$ and $S^\infty$ the unit sphere in the real space $\ell^2$.\ Extrapolating a little bit the concepts found in \cite{schoen}, but still keeping the setting of the general theory developed in \cite{berg}, we will say that a kernel $K: S^m \times S^M \to \mathbb{R}$ is {\em positive definite} if
$$
\sum_{\mu,\nu=1}^n c_\mu c_\nu K((x_\mu, w_\mu),(x_\nu, w_\nu)) \geq 0,
$$
for $n\geq 1$, distinct points $(x_1,w_1), (x_2,w_2), \ldots, (x_n,w_n)$ on $S^m \times S^M$, and real scalars $c_1, c_2, \ldots,$ $c_n$.\ It is {\em strictly positive definite} if it is positive definite and the previous inequalities are strict whenever at least one of the $c_\mu$ is nonzero.\ It is {\em $DC$-strictly positive} if it is positive definite and the previous inequalities are strict whenever the $x$ and the $w$ components of the points are pairwise distinct and at least one of the $c_\mu$ is nonzero.\ Clearly, a strictly positive definite kernel is DC-strictly positive definite but not conversely.\ The symbol $DC$ refers to ``distinct components".

A kernel $K$ acting on $S^m \times S^M$ is {\em isotropic (radial)} if
$$
K((x,z),(y,w))=K_r(x\cdot y,z \cdot w), \quad x,y \in S^{m},\quad z,w \in S^{M},
$$
for some real function $K_r$ on $[-1,1]^2$, where $\cdot$ stands for the inner product of both $\mathbb{R}^{m+1}$ and $\mathbb{R}^{M+1}$.\ In other words, the isotropy of $K$ corresponds to
the property
$$
K((Ax,Bz),(Ay,Bw))=K((x,z),(y,w)), \quad x,y \in S^m, \quad z,w \in S^M, \quad A \in \mathcal{O}^m,\quad B \in \mathcal{O}^M,$$
in which $\mathcal{O}^m$ denotes the orthogonal group in $\mathbb{R}^{m+1}$.\ The kernel $K_r$ in the definition above will be referred to as the {\em isotropic part} of $K$.

Finally, in order to speak of continuity of a kernel as above, we need to assume that all the spheres involved are endowed with their geodesic distances.\ A nice discussion on continuous, isotropic and positive definite kernels on a single sphere, including applications, is made on the recent paper \cite{gneiting}.\ Additional information can be found in references therein.\ As for positive definiteness on a product os spheres, we have found no relevant references to quote, except \cite{jean}.

For $m,M <\infty$, a result proved in \cite{jean} reveals that a real, continuous and isotropic kernel $K$ on $S^m \times S^M$ is positive definite if, and only if, its isotropic part has a double Fourier series representation in the form
$$
K_r(t,s) = \sum_{k,l=0}^{\infty}a_{k,l}P_{k}^{m}(t)P_{l}^{M}(s), \quad t,s \in [-1,1],
$$
in which $\hat{f}_{k,l}\geq 0$, $k,l\in \mathbb{Z}_+$, $P_k^m$ is the Gegenbauer polynomial of degree $k$ with respect to the real number $(m-1)/2$, and
$$
\sum_{k,l=0}^{\infty}a_{k,l}P_{k}^{m}(1)P_{l}^{M}(1)<\infty.
$$
Obviously, this theorem extends a famous result of I. J. Schoenberg (\cite{schoen}) to products of spheres.\ The Gegenbauer polynomials appearing above are discussed in \cite{dai,szego}.\ In particular, one may find there the orthogonality relation for them
$$
\int_{-1}^1 P_n^{m}(t) P_k^{m}(t) (1-t^2)^{(m-2)/2}dt=\frac{\tau_{m+1}}{\tau_{m}}\frac{m-1}{2n+m-1}P_n^{m}(1)\delta_{n,k},
$$
in which $\tau_{m+1}$ is the surface area of $S^{m}$, that is,
$$
\tau_{m+1}:=\frac{2\pi^{(m+1)/2}}{\Gamma((m+1)/2)}.
$$
If one or both spheres coincide with $S^\infty$, the representation theorem described above still holds.\ Indeed, it suffices to replace each Gegenbauer polynomial with the standard monomial of equal degree in the appropriate spots in the expansion for $K_r$.

The results in this paper will apply to the case $m,M\geq 2$ only; in the other cases, some of the corresponding questions are still open while at least one has been settled already (see either comments ahead or \cite{guella}).\ Thus, throughout the paper, we will assume that $m,M\geq 2$.\ For a real, continuous, isotropic and positive definite kernel $K$ on a product of spheres, we can define the set
$$
J_K:= \{(k,l) \in \mathbb{Z}^2: a_{k,l}>0\}
$$
attached to its isotropic part $K_r$ (in the cases we need to work with more then one product simultaneously, we will emphasize the dependance on the dimensions involved by writing $a_{k,l}^{m,M}$ and $J_K(m,M)$).\ One of the results to be proved in this paper is this one.

\begin{thm} \label{main1} Let $K$ be a real, continuous, isotropic and positive definite kernel on $S^m \times S^M$.\ It is $DC$-strictly positive definite if and only if
the sets $\{(k,l) \in J_K: k+l \in 2\mathbb{Z}_+\}$ and $\{(k,l) \in J_K: k+l \in 1+2\mathbb{Z}_+\}$ are infinite.
\end{thm}

For a continuous, isotropic and positive definite kernel $K$, we can also also define the sets
\begin{align*}
J_K^{0,0}& = J_K \cap (2\mathbb{Z}_+ \times 2\mathbb{Z}_+) \\
J_K^{0,1}& = J_K\cap [2\mathbb{Z}_+  \times (1+2\mathbb{Z}_+)] \\
J_K^{1,0}& = J_K\cap [(1+2\mathbb{Z}_+) \times 2\mathbb{Z}_+] \\
J_K^{1,1}& = J_K \cap [(1+2\mathbb{Z}_+) \times (1+2\mathbb{Z}_+)].
\end{align*}
so that
$$
J_K=J_K^{0,0} \cupdot J_K^{0,1} \cupdot J_K^{1,0} \cupdot J_K^{1,1}.$$
The second main contribution in this paper is this one.

\begin{thm} \label{main2} Let $K$ be a real, continuous, isotropic and positive definite kernel on $S^m \times S^M$.\ It is strictly positive definite if and only if for each pair $(i,j)$, there exists a sequence $\{(k_r^{i,j},l_r^{i,j})\}_{r \in \mathbb{Z}_+}$ in $J_K^{i,j}$ so that $\lim_{r \to \infty}k_r^{i,j}=\lim_{r\to \infty}l_r^{i,j}=\infty$.
\end{thm}

A version of this theorem in the case in which $m=M=1$ was obtained in \cite{guella}.

The paper proceeds as follows.\ In Section 2, we present an alternative formulation for the concept of strict positive definiteness of a real, continuous, isotropic and positive definite kernel on $S^m \times S^M$.\ The proof of Theorem \ref{main1} is a consequence of the results proved in Section 3, while Theorem \ref{main2} follows from the results to be presented in Section 4.

\section{Strict positive definiteness: technical results}

In order to verify that a continuous, isotropic and positive definite kernel $K$ on $S^{m}\times S^{M}$ is either strictly positive definite or $DC$-strictly positive definite on $S^m \times S^M$, one needs to deal with the positive definiteness of matrices $A$ having $\mu\nu$-entries in the form
$$
A_{\mu\nu}=K_r(x_\mu \cdot x_{\nu},w_\mu \cdot w_\nu),
$$
for $n$ distinct points $(x_1, w_1), (x_2, w_2), \ldots, (x_n, w_n)$ on $S^{m} \times S^M$.\ In other words, no matter what the data above is, with or without the additional requirement in the definition of $DC$-strict positive definiteness, one needs to conclude that $c=(c_1,c_2,\ldots,c_n)=0$ from the equality
$$
c^t A c=\sum_{\mu=1}^n \sum_{\nu=1}^n c_{\mu}c_{\nu}K_r(x_\mu \cdot x_{\nu},w_\mu \cdot w_\nu)=0.
$$
In particular, one needs to extract one vector equality ($n$ scalar equalities) from just one scalar equation.\ Depending on the situation, this may be an indigestible task, perhaps impossible.\ Roughly speaking, the results in this section will provide a concise option to change the equation $c^tAc=0$ into a set of functional equations, a procedure that allows the task mentioned above to be performed in a quite more reasonable way.\ The notation $c^tAc$ will always refer to the setting explained above.

We will write $\mathcal{H}_k^m$ to denote the space of all spherical harmonics of degree $k$ in $m+1$ variables with dimension $d(k,m)$ and $\{Y_{k,1}^m, Y_{k,2}^m,\ldots,Y_{k,d(k,m)}^m\}$ to denote a basis for that space.\ In particular, $d(k,m)$ stands for the dimension of the space $\mathcal{H}_k^m$.\ The addition theorem asserts that (\cite[p. 10]{dai})
$$
P_k^{m}(x\cdot y)=\frac{m-1}{2k+m-1}\sum_{j=1}^{d(k,m)}Y_{k,j}^m(x)Y_{k,j}^m(y), \quad x,y \in S^{m}.
$$
More information on spherical harmonics can be found in \cite{atkinson,groemer, muller}.

\begin{prop} \label{tec1} Let $K$ be a nonzero, real, continuous, isotropic and positive definite kernel on $S^m \times S^M$.\ For distinct points $(x_1, w_1), (x_2, w_2), \ldots, (x_n, w_n)$ on $S^{m} \times S^M$ and real numbers $c_1, c_2,\ldots, c_n$, the following assertions are equivalent:\\
$(i)$ $c^{t}Ac=0$;\\
$(ii)$ The equality
$$
\sum_{\mu=1}^{n}c_{\mu}P_{k}^{m}(x_{\mu}\cdot x)P_{l}^{M}(w_{\mu}\cdot w)=0
$$
holds for $(k,l) \in J_K$, $x\in S^{m}$ and $w \in S^{M}$.
\end{prop}
\pf The addition theorem and the representation for the isotropic part $K_r$ of $K$ justify the equality
$$
c^{t}Ac=\sum_{k,l=0}^{\infty}\frac{m-1}{2k+m-1}\frac{M-1}{2l+M-1}a_{k,l}\sum_{i=1}^{d(k,m)}\sum_{j=1}^{d(l,M)}\left|\sum_{\mu=1}^n c_\mu Y_{k,i}^m(x_\mu)Y_{l,j}^M(w_\mu)\right|^2.
$$
Hence, $c^{t}Ac=0$, if and only if,
$$
\sum_{\mu=1}^n c_\mu Y_{k,i}^m(x_\mu)Y_{l,j}^M(w_\mu)=0, \quad (k,l) \in J_K, \quad i=1,2,\ldots,d(k,m), \quad j=1,2,\ldots,d(l,M).
$$
Multiplying the previous equality by $Y_{k,i}^m(x)Y_{l,j}^M(w)$, adding up on $i$ and $j$ and using the addition theorem once again leads to the statement in $(ii)$.\ Conversely, if $(ii)$ holds, the addition formula implies that
$$
\sum_{j=1}^{d(l,M)}\left[\sum_{i=1}^{d(k,m)}\sum_{\mu=1}^n c_\mu Y_{k,i}^m(x_\mu)Y_{k,i}^m(x)Y_{l,j}^M(w_\mu)\right]Y_{l,j}^M=0, \quad x \in S^{m}, \quad (k,l) \in J_K.
$$
Since $\{Y_{l,1}^M, Y_{l,2}^M,\ldots,Y_{l,d(l,M)}^M\}$ is a basis of $\mathcal{H}_l^M$ for all $l$, then
$$
\sum_{i=1}^{d(k,m)}\left[\sum_{\mu=1}^n c_\mu Y_{k,i}^m(x_\mu)Y_{l,j}^M(w_\mu)\right]Y_{k,i}^m=0, \quad j=1,2,\ldots,d(l,M),\quad  (k,l) \in J_K.
$$
Likewise, since $\{Y_{k,1}^m, Y_{k,2}^m,\ldots,Y_{k,d(k,m)}^m\}$ is a basis of $\mathcal{H}_k^m$ for all $k$, and taking into account the first equality in the proof, the equality
$c^tAc=0$ follows.\eop\vspace*{3mm}

The theorems below are obvious consequences of the proposition.

\begin{thm}
Let $K$ be a real, continuous, isotropic and positive definite kernel on $S^m \times S^M$.\ The following assertions are equivalent:\\
$(i)$ $K$ is $DC$-strictly positive definite;\\
$(ii)$ If $n\geq 1$ and $(x_1, w_1), (x_2, w_2), \ldots, (x_n, w_n)$ are distinct points on $S^{m} \times S^M$ possessing distinct components, then the only solution of the system
$$
\sum_{\mu=1}^{n}c_{\mu}P_{k}^{m}(x_{\mu}\cdot x)P_{l}^{M}(w_{\mu}\cdot w)=0, \quad (k,l) \in J_K,\quad x\in S^{m}, \quad w \in S^{M},$$
is $c=0$.
\end{thm}

\begin{thm} \label{caracSPD}
Let $K$ be a real, continuous, isotropic and positive definite kernel on $S^m \times S^M$.\ The following assertions are equivalent:\\
$(i)$ $K$ is strictly positive definite;\\
$(ii)$ If $n\geq 1$ and $(x_1, w_1), (x_2, w_2), \ldots, (x_n, w_n)$ are distinct points on $S^{m} \times S^M$, then the only solution of the system
$$
\sum_{\mu=1}^{n}c_{\mu}P_{k}^{m}(x_{\mu}\cdot x)P_{l}^{M}(w_{\mu}\cdot w)=0, \quad (k,l) \in J_K,\quad x\in S^{m}, \quad w \in S^{M},$$
is $c=0$.
\end{thm}

\section{$DC$-strict positive definiteness on $S^m \times S^M$, $m,M >1$}

The results to be proved in this section will justify the first main contribution of the paper, namely, Theorem \ref{main1}.\ If one compares the concept of $DC$-strict positive definiteness introduced here with that one of strict positive definiteness on a single sphere originally considered in \cite{cheney}, it becomes quite clear that the first one is a true two-dimensional version of the second one.\ As so, one should expect that  a characterization of $DC$-strict positive definiteness should be the natural extension of that for strict positive definiteness on just one sphere in \cite{chen}.\ Thus, to a certain extension, the results here are an upgrade of those in \cite{chen}.\ In contrast with the proofs in \cite{chen}, the arguments to be used here make no use of spherical coordinate systems in either sphere involved.

The necessity of the condition quoted in Theorem \ref{main1} is easily verified following basic estimations on the rank of the matrices $A$ mentioned in the previous section, for special point distributions on $S^m \times S^M$.\ In order to explain that, we will use the fact that the Gegenbauer polynomials of even degree are even functions while Gegenbauer polynomials of odd degree are odd functions.\ We include the proof in the case in which both $m$ and $M$ are finite; the reader can verify himself that it works on the other cases too, after making the obvious modifications.

\begin{prop}\label{preproj}
Let $K$ be a real, nonzero, continuous, isotropic and positive definite kernel on $S^m \times S^M$.\ In order that it be $DC$-strictly positive definite it is necessary that $\{k+l: (k,l) \in J_K\}$ contains infinitely many even and infinitely many odd integers.
\end{prop}
\pf Let us assume that $\{(k,l) \in J_K: k+l \in 2\mathbb{Z}_+\}$ is nonempty and finite.\ Hence, we can write
$$
K_r(t,s)=\sum_{k=0}^{N}\sum_{l=0}^{N}a_{k,l}P_{k}^{m}(t)P_{l}^{M}(s)+\sum_{(k,l)\in K_1}a_{k,l}P_{k}^{m}(t)P_{l}^{M}(s),\quad t,s \in [0,1],
$$
in which $N$ is a nonnegative integer and every element $(k,l)$ of $K_1$ satisfies $k+l \in 1+2\mathbb{Z}_+$.\ Now, pick $2n$ distinct points $x_1,x_2, \ldots, x_{2n}$ in $S^{m}$ and $2n$ distinct points $w_1,w_2, \ldots, w_{2n}$ in $S^{M}$, chosen so that $x_{n+j}=-x_j$ and $w_{n+j}=-w_j$, for $j=1,2, \ldots,n$.\  The $2n \times 2n$ matrix $A=[K_r(x_\mu \cdot x_\nu, w_\mu \cdot w_\nu)]$ breaks down into a sum $B+C$, in accordance with the decomposition of $K_r$ introduced above.\ If we define $c_j$ as the vector having its $j$-th and $(n+j)$-th components equal to 1, then it is straightforward to verify that $\{c_1, c_2,\ldots,c_n\}$ is a subset of the kernel of $C$.\ In particular, the rank of $C$ is at most $n$.\ On the other hand, since the $2n \times 2n$ matrix $[x_\mu \cdot x_\nu]$ has rank at most $m+1$, its Schur product $[(x_\mu \cdot x_\nu)^j]$ has rank at most $(m+1)^j$.\ Since $P_k^{m}$ is a polynomial of degree $k$, the matrix
$[P_k^m(x_\mu \cdot x_\nu)]$ has rank at most $1+m+1+\cdots+(m+1)^k\leq (N+1)(m+1)^N$, whenever $k\leq N$.\ It is now clear that the rank of $B$ is at most $(N+1)^4[(m+1)(M+1)]^N$.\ Thus, if $n> (N+1)^4[(m+1)(M+1)]^N$, then $A$ is not of full rank $2n$ and $K$ cannot be $DC$-strictly positive definite on $S^{m} \times S^{M}$.\ If the assumption is the nonemptiness and finiteness of $\{(k,l) \in J_K: k+l \in 1+2\mathbb{Z}_+\}$, then a similar procedure can be applied.\ Finally, if one of the sets $\{(k,l) \in J_K: k+l \in 2\mathbb{Z}_+\}$ and $\{(k,l) \in J_K: k+l \in 1+2\mathbb{Z}_+\}$ is empty, the procedure above works using an arbitrary positive integer $n$ and the the same choices of points. \eop\vspace*{3mm}

An alternative proof for Proposition \ref{preproj} can be achieved via Proposition 4.4 in \cite{jean}.\ Indeed, if $K$ is a real, continuous, isotropic and $DC$-strictly positive definite kernel on $S^m \times S^M$, then that proposition implies that $t \in [-1,1] \to K_r(t,t)$ is a continuous, isotropic and strictly positive definite kernel on $S^{m\wedge M}$, in which $m\wedge M:=\min\{m,M\}$.\ Since the linearization formula
$$P_k^m(t) P_l^M(t)=\sum_{j=0}^{\lfloor (k+l)/2\rfloor} \alpha_j^{m,M}(k,l) P_{k+l-2j}^{m\wedge M}(t), \quad \alpha_j^{m,M}(k,l)>0,$$
 holds for all $k$ and $l$, it is promptly seen that
$$K_r(t,t)=\left[\sum_{k+l\in 1+ 2\mathbb{Z}_+}+\sum_{k+l\in 2\mathbb{Z}_+}\right]a_{k,l}\sum_{j=0}^{\lfloor (k+l)/2\rfloor} \alpha_j^{m,M}(k,l)P_{k+l-2j}^{m\wedge M} (t)\quad t \in [-1,1].$$
Invoking the main theorem in \cite{chen}, we now see that
$$a_{k,l}\sum_{j=0}^{\lfloor (k+l)/2\rfloor} \alpha_j^{m,M}(k,l) >0$$
for infinitely many pairs $(k,l)$ with $k+l$ even and infinitely many pairs with $k+l$ odd.\ But that implies the necessary condition in Proposition \ref{preproj}.

Next, we normalize the Gegenbauer polynomials by writing
$$
R_k^{m}=\frac{P_k^{m}}{P_k^{m}(1)}, \quad k\in \mathbb{Z}_+.
$$
Since $|P_k^{m}(t)|\leq P_k^{m}(1)$, $t \in [-1,1]$, we immediately have that
$$|R_k^m(t)| \leq 1, \quad k \in\mathbb{Z}_+, \quad t \in [-1,1].$$
Two more specific properties of the Gegenbauer polynomials are included in the lemma below.

\begin{lem} \label{gege} The following assertions hold:\\
$(i)$ $|R_k^m(t)| =1$ if and only if either $t=1$ or $t=-1$.\\
$(ii)$ $ \lim_{k \to \infty} R_{k}^{m}(t)=0$, $t \in (-1,1)$.
\end{lem}
\pf Assertion $(i)$ can be deduced from the fact that Gegenbauer polynomials can be expanded in terms of Tchebyshev polynomials (\cite[p.59]{askey}): if $k$ is a nonnegative integer, then there
are positive constants $c_k^m(j)$, $j=0,1,\ldots,k$, such that
$$P_k^m(t)=\sum_{j=0}^{\lfloor k/2\rfloor} c_k^m(j) P_{k-2j}^0(t), \quad t \in [-1,1].$$
As for $(ii)$,
it is implied by a quite general inequality for Jacobi polynomials described in \cite[p.196]{szego} (see also \cite[p.416]{dai}). \eop

Lemma \ref{gege}-$(ii)$ plays an important role in this paper.\ Since it does not hold in the case $m=1$, the characterization for strict positive definiteness and $DC$ strict positive definiteness in the cases when at least one of the spheres is a circle needs to follow a different pattern.\ This is one of the reasons why those cases cannot be included in this paper.

We are about ready to prove the converse of the previous proposition.

\begin{thm} \label{mainjean}
Let $K$ be a real, continuous, isotropic and positive definite kernel on $S^m \times S^M$.\ In order that it be $DC$-strictly positive definite it is sufficient that $\{k+l: (k,l) \in J_K\}$ contains infinitely many even and infinitely many odd integers.
\end{thm}
\pf Assume $\{(k,l) \in J_K: k+l \in 2\mathbb{Z}_+\}$ and $\{(k,l) \in J_K: k+l \in 1+2\mathbb{Z}_+ \}$ are infinite, let $n$ be a positive integer, $x_1,x_2, \ldots, x_n$ distinct points in $S^M$ and $w_1, w_2, \ldots, w_n$ distinct points on $S^M$.\ As before, write $A$ to denote the $n \times n$ matrix with entries $A_{\mu\nu}=K_r(x_\mu \cdot x_\nu, w_\mu \cdot w_\nu)$.\
We will show that the equality $c^tAc=0$ implies $c=0$.\ Applying Proposition \ref{tec1}, we know that the equality $c^tAc=0$ corresponds to
$$
\sum_{\mu=1}^{n}c_{\mu}P_{k}^{m}(x_{\mu}\cdot x)P_l^{M}(w_{\mu} \cdot w)=0,\quad (k,l) \in J_K, \quad x\in S^{m}, \quad w\in S^{M}.
$$
Next, for $\alpha$ in $\{1,2,\ldots,n\}$ fixed, we will conclude that $c_\alpha=0$ via specific choices for
the points $x \in S^{m}$ and $w \in S^{M}$ in the equality above and with the help of some special computations.\ There are 5 distinct cases to be considered.\\
\underline{Case $1$} {\em No $x_\mu$ is antipodal of $x_{\alpha}$ and no $w_\mu$ is antipodal of $w_{\alpha}$.}\\
Here, the choices are $x=x_{\alpha}$ and $w=w_{\alpha}$, the resulting equation being
$$
c_{\alpha} P_k^{m}(1) P_l^{M}(1)+\sum_{\mu \neq \alpha} c_\mu  P_{k}^{m}(x_{\mu}\cdot
x_{\alpha})P_l^{M}(w_{\mu} \cdot w_{\alpha})=0, \quad (k,l) \in J_K,
$$
that is,
$$
c_{\alpha}+\sum_{\mu \neq \alpha} c_\mu  R_{k}^{m}(x_{\mu}\cdot
x_{\alpha})R_l^{M}(w_{\mu} \cdot w_{\alpha})=0, \quad (k,l) \in J_K.
$$
Due to our assumption on $\{k+l: (k,l) \in J_K\}$, we can select a sequence $\{(k_r,l_r)\}\subset J_K$ so that either $\lim_{r \to \infty}k_r=\infty$ or $\lim_{r\to \infty} l_r=\infty$.\ Introducing the sequence in the previous equality, observing that
$$x_\mu \cdot x_\alpha \neq \pm 1 \neq w_{\mu}\cdot w_{\alpha},\quad \mu \neq \alpha $$ and calculating the limit as $r \to \infty$ with a help of Lemma \ref{gege}-$(ii)$, we obtain
$$
0 =c_\alpha +\lim_{r \to \infty} \sum_{\mu \neq \alpha} c_\mu  R_{k_{r}}^{m}(x_{\mu}\cdot
x_{\alpha})R_{l_{r}}^{M}(w_{\mu}\cdot w_{\alpha})=c_\alpha.
$$
\underline{Case $2$} {\em No $x_\mu$ is antipodal of $x_{\alpha}$ and some $w_{\beta}$ is antipodal of $w_{\alpha}$}.\\
The same choice used in the previous case leads to
$$
c_{\alpha}  +  (-1)^l c_{\beta} R_k^{m}(x_{\beta} \cdot x_{\alpha})+\sum_{\mu\neq \alpha,\beta}c_\mu  R_{k}^{m}(x_{\mu}\cdot
x_{\alpha})R_l^{M}(w_{\mu}\cdot w_{\alpha})=0, \quad (k,l) \in J_K.$$
If there exists a sequence $\{(k_r,l_r)\}\subset J_K$ so that $ \lim_{r \to \infty} k_r=\infty$, then Lemma \ref{gege}-$(ii)$ implies that
$$
0=c_\alpha+\lim_{r\to \infty} (-1)^{l_r} c_{\beta}R_{k_r}^{m}(x_{\beta} \cdot x_{\alpha})+
\lim_{r \to \infty}\sum_{\mu\neq \alpha,\beta}c_\mu  R_{k_r}^{m}(x_{\mu}\cdot
x_{\alpha})R_{l_r}^{M}(w_{\mu}\cdot w_{\alpha})=c_\alpha.
$$
Otherwise, we can select a nonnegative integer $k$ and a  sequence  $\{l_r\} \subset \mathbb{Z}_+$ fulfilling the following requirements: the parity among the elements in the sequence is the same, $\{(k, l_r) \} \subset J_K$, and $\lim_{r \to \infty}l_r= \infty$.\ Inserting this sequence in the initial equation of the case, we obtain
$$
c_{\alpha}  +  (-1)^{l_{r}} c_{\beta} R_k^{m}(x_{\beta} \cdot x_{\alpha})+\sum_{\mu\neq \alpha,\beta}c_\mu  R_{k}^{m}(x_{\mu}\cdot
x_{\alpha})R_{l_{r}}^{M}(w_{\mu}\cdot w_{\alpha})=0.
$$
Letting $r\to \infty$, we conclude that
$$
c_\alpha +c_\beta (-1)^pR_{k}^{m}(x_\beta \cdot x_\alpha)=0,
$$
in which $p$ is either even or odd, depending on the parity of the elements of the sequence $\{l_r\}$.\ To conclude the arguments, we consider a second choice for $x$ and $w$, namely, $x=x_{\beta}$ and $w=w_{\alpha}$ and insert the same sequence to reach
$$
c_{\alpha} R_k^{m}(x_\alpha \cdot x_\beta) +  (-1)^{l_{r}} c_{\beta} +\sum_{\mu\neq \alpha,\beta}c_\mu  R_{k}^{m}(x_{\mu}\cdot
x_{\beta})R_{l_{r}}^{M}(w_{\mu}\cdot w_{\alpha})=0.
$$
Letting $r \to \infty$ once again, now leads to
$$
c_{\alpha} R_{k}^{m}(x_\alpha \cdot x_\beta) +  (-1)^p c_{\beta} =0.
$$
The two relations between $c_\alpha$ and $c_\beta$ deduced above produce the single equation
$$
c_\alpha\left[1- R_{k}^{m}(x_\alpha \cdot x_\beta)^{2} \right]=0.
$$
Taking into account that $x_{\beta }$ and $x_{\alpha}$ are distinct and not antipodal, Lemma \ref{gege}-$(i)$ implies that $c_\alpha=0$.\\
\underline{Case 3} {\em Some $x_{\beta}$ is antipodal of $x_{\alpha}$ and no $w_{\mu}$ is antipodal of $w_{\alpha}$}.\\ This case is similar to the previous one.\\
\underline{Case 4} {\em Some $x_{\beta}$ is antipodal of $x_{\alpha}$ and $w_{\beta}$ is antipodal of $w_{\alpha}$}.\\
The choices $x=x_{\alpha}$ and $w = w_{\alpha}$ leads to
$$
c_{\alpha}  +  (-1)^{k+l} c_{\beta}+\sum_{\mu\neq \alpha,\beta}c_\mu  R_{k}^{m}(x_{\mu}\cdot
x_{\alpha})R_l^{M}(w_{\mu}\cdot w_{\alpha})=0, \quad (k,l)\in J_K.
$$
Our assumption on $J_K$ allows the selection of two subsequences $\{(k_r,l_r)\}$ and $\{(k_s,l_s)\}$ of $J_K$ with $k_r+l_r\in 2\mathbb{Z}_+$ for all $(k,l)$, $k_s+l_s \in 1+2\mathbb{Z}_+$ for all $(k,l)$, $\lim_{r \to \infty}k_r+ l_r =\infty$, and $\lim_{s \to \infty}k_s+ l_s =\infty$.\ Introducing these two sequences in the previous equation and letting $r,s \to \infty$ lead to
$$c_{\alpha}+c_{\beta}=c_{\alpha}-c_{\beta}=0.$$ In particular, $c_{\alpha}=0$.\\
\underline{Case 5} {\em Some $x_{\beta_{1}}$ is antipodal of $x_{\alpha}$, some $w_{\beta_{2}}$ is antipodal of $w_{\alpha}$, and $\beta_{1} \neq  \beta_{2}$. }\\
This case requires the consideration of subcases.\ If there exists a sequence $\{(k_r,l_r)\}\subset J_K$ for which $ \lim_{r \to \infty} k_r= \lim_{r \to \infty} l_r=\infty$, then choosing $x=x_{\alpha}$ and $w = w_{\alpha}$ in the original equality and inserting the sequence leads to
$$
c_{\alpha}  +  (-1)^{k_{r}} c_{\beta_{1}} R_{l_r}^{m}(w_{\beta_{1}} \cdot w_{\alpha})+ (-1)^{l_{r}} c_{\beta_{2}} R_{k_r}^{m}(x_{\beta_{2}} \cdot x_{\alpha})+\sum_{\mu\neq \alpha,\beta_{1}, \beta_{2}}c_\mu  R_{k_r}^{m}(x_{\mu}\cdot
x_{\alpha})R_{l_{r}}^{M}(w_{\mu}\cdot w_{\alpha})=0.
$$
Letting $r \to \infty$ we reach $c_\alpha=0$.\ If there exists no such sequence in $J_K$, our assumptions on $J_K$ leaves us with two possibilities: either there is a fixed $k$ and a sequence $\{l_{r}\}$ so that all the $l_r$ have the same parity $p$, $\{(k,l_{r})\} \subset J_K$ and $\lim_{r \to \infty }l_{r} = \infty$ or there is a fixed $l$ and a sequence $\{k_{r}\}$ so that all the $k_r$ have the same parity $p$, $\{(k,l_{r})\} \subset J_K$ and $\lim_{r \to \infty }k_{r} = \infty$.\ We will proceed with the first possibility, the other one being similar.\ The choice $x=x_{\alpha}$ and $w = w_{\alpha}$ in the original equation and the insertion of the sequence of the subcase provides the relation
$$
c_{\alpha} +  (-1)^k c_{\beta_{1}}R_{l_{r}}^{M}(w_{\beta_{1}}\cdot w_{\alpha})  +   (-1)^{l{_{r}}} c_{\beta_{2}}R_{k}^{m}(x_{\beta_{2}}\cdot x_{\alpha})
+ \sum_{\mu\neq \alpha,\beta_{1}, \beta_{2}}c_\mu  R_{k}^{m}(x_{\mu}\cdot x_{\alpha})R_{l_{r}}^{M}(w_{\mu}\cdot w_{\alpha})=0,
$$
while the choice $x=x_{\beta_{2}}$ and $w=w_{\alpha}$ yields
$$
c_{\alpha}R_{k}^{m}(x_{\beta_{2}}\cdot x_{\alpha})  +  c_{\beta_{1}}R_{k}^{m}(x_{\beta_{1}}\cdot x_{\beta_{2}})R_{l_{r}}^{M}(w_{\beta_{1}}\cdot w_{\alpha})
  + (-1)^{l{_r}} c_{\beta_{2}} + \sum_{\mu\neq \alpha,\beta_{1}, \beta_{2}}c_\mu  R_{k}^{m}(x_{\mu}\cdot x_{\beta_{2}})R_{l_{r}}^{M}(w_{\mu}\cdot w_{\alpha})=0.
$$
Letting $r \to \infty$ in both equations, we deduce that
$$
c_\alpha +c_{\beta_{2}} (-1)^p R_{k}^{m}(x_{\beta_{2}} \cdot x_{\alpha} )=c_{\alpha}R_{k}^{m}(x_{\beta_{2}}\cdot x_{\alpha}) + (-1)^{p} c_{\beta_{2}}=0.
$$
It is now clear that
$$
c_\alpha\left[1- R_{k}^{m}(x_\alpha \cdot x_{\beta_{2}})^{2} \right]=0.
$$
Since $x_{\beta_{2} }$ and $x_{\alpha}$ are distinct and not antipodal, Lemma \ref{gege}-$(ii)$ implies that $c_\alpha=0$.\eop

The technique presented in proof of Theorem \ref{mainjean} can be used in the writing of an alternative proof of the main theorem in \cite{chen}, one that does not depend upon special coordinate systems on the sphere.

Another important fact to be noticed at this time is that the proof of Theorem \ref{mainjean} cannot be modified in order to produce a sufficient condition for the plain strict positive definiteness of the kernel $K$.\ Indeed,
if we relax our assumptions and permit repetitions among either the $x_\mu$ or the $w_\mu$, then we will no longer be able to use Lemma \ref{gege}-$(ii)$.\ Therefore, a different approach needs to be implemented in the general case, as the reader shall see in the next section.

Lemma \ref{gege} does not hold when $m=1$.\ Therefore, the cases in which either $m=1$ or $M=1$ cannot be included in the statement of Theorem \ref{mainjean}.\ A characterization for $DC$-strict positive definiteness in the cases in which at least one of the spheres is a circle remains open at this time.

To close the section we will ratify Theorem \ref{main1} in the case when either $m$ or $M$ is infinite.\ We need two technical lemmas, the first one originally proved in \cite{bingham}.

\begin{lem} \label{decompo} Let $m$ be an integer at least 2.\ If $k$ is a nonnegative integer, then there are positive integers
$c(k,m,j)$, $j=0,1, \ldots, \lfloor k/2\rfloor$, so that
$$
t^{k}=\sum_{j=0}^{\lfloor k/2\rfloor} c(k,m,j)P^{m}_{k-2j}(t), \quad t\in [-1,1].
$$
\end{lem}

As for the other one, we sketch the proof.

\begin{lem} \label{spdfinite} Let $K$ be a continuous, isotropic and positive definite kernel on $S^\infty \times S^{M}$.\ The following assertions are equivalent:\\
$(i)$ $K$ is $DC$-strictly positive definite on $S^\infty \times S^{M}$;\\
$(ii)$ $K$ is $DC$-strictly positive definite on $S^{m} \times S^{M}$, $m=2,3,\ldots$.
\end{lem}\label{equivinfty}
\pf One implication uses the fact that each $S^m$ is isometrically embedded in $S^\infty$.\ Hence, so does $S^m \times S^M$ in $S^\infty \times S^M$.\ As for the other, one needs to use the fact that the linear span of a subset of $S^\infty$ containing $n$ points can be seen as a subset of a copy of $S^n$ isometrically embedded in $S^\infty$. \eop\vspace*{3mm}

The isotropic part $K_r$ of a real, continuous, isotropic and positive definite kernel $K$ on $S^m \times S^M$ is also the isotropic part of a continuous, isotropic and positive definite kernel on $S^n \times S^M$, for $n\leq m$.\ As so, if the dimensions are to be emphasized, we will write
$K_r^{m,M}$ instead of $K_r$, $a_{k,l}^{m,M}$ for the Fourier coefficients of $K_r^{m,M}$ and $J_K(m,M)$ for the index set $J_K$.

Here is the version of Theorem \ref{main1} when $m=\infty$.

\begin{thm} \label{mainpdinfty} Let $K$ be a continuous, isotropic and positive definite kernel on $S^\infty \times S^M$.\ In order that it be $DC$-strictly positive definite it is necessary and sufficient that both sets $\{(k,l) \in J_K(\infty,M): k+l \in 2\mathbb{Z}_+\}$ and $\{(k,l) \in J_K(\infty,M): k+l \in 1+2\mathbb{Z}_+ \}$ be infinite.
\end{thm}
\pf If $m$ is a positive integer at least 2, we can use Lemma \ref{decompo} to write
$$
K_r^{\infty,M}(t,s)=\sum_{k=0}^{\infty} \sum_{l=0}^{\infty}\left[a_{k,l}^{\infty,M}\sum_{j=0}^{\lfloor k/2\rfloor} c(k,m,j) \right]P^{m}_{k-2j}(t) P_{l}^{M}(s), \quad t,s \in [-1,1].
$$
Obviously, the formula can be re-written as
$$
K_r^{\infty,M}(t,s)=\sum_{k=0}^{\infty} \sum_{l=0}^{\infty}\left[\sum_{n=0}^{\infty} c(k + 2n,m,n)a_{k+2n,l}^{\infty,M}\right]P^{m}_{k}(t) P_{l}^{M}(s), \quad t,s \in [-1,1].
$$
Due to the orthogonality relation for the Gegenbauer polynomials, it is not hard to see that the functions $(t,s) \in [-1,1]^2 \to P^{m}_{k}(t) P_{l}^{M}(s)$ form an orthogonal system in the space $L^2([-1,1],w_{m,M})$, in which
$$
w_{m,M}(t,s)=(1-t^2)^{(m-2)/2}(1-s^2)^{(M-2)/2},\quad t,s \in [-1,1].
$$
Hence, the expansion for $K_r^{\infty,M}$ above corresponds to the series representation of $K_r^{m,M}$.\ In particular,
the number $c_{k,l}^{m,M}:=\sum_{n=0}^{\infty} c(k + 2n,m,n)a_{k+2n,l}^{\infty,M}$ is a positive multiple of $a^{m,M}_{k,l}$.\ Now, we have a solid terrain on which the proof of the theorem itself can be drafted.\ If either set in the statement of the theorem is finite, then we can pick $m\geq 2$ of our choice and conclude that either $\{(k,l) \in J_K(m,M): k+l \in 2\mathbb{Z}_+\}$ or $\{(k,l) \in J_K(m,M): k+l \in 1+2\mathbb{Z}_+ \}$ is finite.\ That being the case, Theorem \ref{main1} reveals that $K$ is not $DC$-strictly positive definite on $S^m \times S^M$.\ Due to Lemma \ref{spdfinite}, we now see that $K$ cannot be $DC$-strictly positive definite on $S^\infty \times S^M$.\ Conversely, if both sets in the statement of the theorem are infinite, then no matter what $m\geq 2$ we pick, we will have that $c_{k,l}^{m,M}>0$, for infinitely many elements in both $\{(k,l) \in J_K(m,M): k+l \in 2\mathbb{Z}_+\}$ and $\{(k,l) \in J_K(m,M): k+l \in 1+2\mathbb{Z}_+ \}$.\ Thus, the same Theorem \ref{mainjean} implies that $K$ is $DC$-strictly positive definite on every $S^m \times S^M$.\ An application of Lemma \ref{spdfinite} reveals that $K$ is strictly positive definite on $S^\infty \times S^{M}$.\eop

Adapting the arguments used in the proof of Theorem \ref{mainpdinfty} and invoking the theorem itself, one can deduce the following additional result.

\begin{thm} \label{mainpdinfty1}
Let $K$ be a continuous, isotropic and positive definite kernel on $S^\infty \times S^\infty$.\ In order that it be $DC$-strictly positive definite it is necessary and sufficient that both sets $\{(k,l) \in J_K(\infty,\infty): k+l \in 2\mathbb{Z}_+\}$ and $\{(k,l) \in J_K(\infty,\infty): k+l \in 1+2\mathbb{Z}_+ \}$ be infinite.
\end{thm}

A closer look at the proof of Theorem \ref{mainpdinfty} reveals that one can deduce a characterization for plain strict positive definiteness on $S^m \times S^M$, in the case when either $m=\infty$ or $M=\infty$, as long as a similar characterization for strict positive definiteness on $S^m \times S^M$, $m,M <\infty$, is available.

\section{Plain strict positive definiteness}

The intended goal in this section is to present a proof for Theorem \ref{main2}.\ While the characterization for $DC$-strict positive definiteness described in Section 3 had some resemblance with that for strict positive definiteness on a single sphere, a characterization for strict positive definiteness on $S^m \times S^M$ demands quite different arguments.

We begin with the notion of antipodal free sets in $S^m \times S^M$.\ A subset $\{(x_{1}, w_{1}),(x_{2}, w_{2}),$ $ \dots , (x_{n}, w_{n})\}$ of $S^{m}\times S^{M}$ is {\em antipodal free} if there are no pairs of antipodal points among the $x_\mu$ and among the $w_\mu$.\ We observe that, even the points in an antipodal free set being distinct, some of the $x$ components and some of the $w$ components may be equal.

The lemma below concerns the extraction of antipodal free sets from a given set.\ Since the result is quite elementary, the proof will be omitted.

\begin{lem} \label{antipodal} If $(x_{1}, w_{1}),(x_{2}, w_{2}), \dots , (x_{n}, w_{n})$ are distinct points on $S^{m}\times S^{M}$, then there exists $p\leq n$ and an antipodal free subset $\{(x'_{1}, w'_{1}),(x'_{2}, w'_{2}), \dots , (x'_{p}, w'_{p})\}$ of $S^{m}\times S^{M}$ so that
$$
\{(x_{\mu}, w_{\mu} ) :\mu=1,2,\ldots,n\} \subset \{ (\pm x'_{1}, \pm w'_{1}),(\pm x'_{2}, \pm w'_{2}), \dots , (\pm x'_{p}, \pm w'_{p})  \}. $$
\end{lem}

Needless to say that the antipodal free set provided by the previous lemma is not unique.

Next, we present an alternative formulation for Theorem \ref{caracSPD}.\ The letter $I$ in the statement of the theorem will refer to the set $\{(0,0), (1,0),(0,1),(1,1)\}$.

\begin{thm} \label{caracJean1} Let $K$ be a real, continuous, isotropic and positive definite kernels on $S^m\times S^M$.\ The following assertions are equivalent:\\
$(i)$ $K$ is strictly positive definite on $S^{m}\times S^{M}$;\\
$(ii)$ If $n$ is a positive integer and $\{(x_{1}, w_{1}),(x_{2}, w_{2}), \dots , (x_{n}, w_{n})\}$ is an antipodal free subset of $S^{m}\times S^{M}$, then the only solution $(c_{1}^{0,0},\ldots,c_{n}^{0,0},c_{1}^{1,0},\ldots,c_{n}^{1,0},
 c_{1}^{0,1}, \ldots, c_{n}^{0,1},c_{1}^{1,1}, \ldots, c_{n}^{1,1})$ to the system of equations
$$
\sum_{\mu =1}^{n} \left[\sum_{(i,j) \in I}(-1)^{ik+jl}c_{\mu}^{i,j}\right] P_{k}^{m}( x_{\mu}\cdot x)P_{l}^{M}( w_{\mu} \cdot w)=0, \quad (x,w) \in S^m \times S^M, \quad (k,l) \in J_K,
$$
is the zero vector.
\end{thm}
\pf Assume $(i)$ holds.\ The set
$$\{(\pm x_{1}, \pm w_{1}),(\pm x_{2}, \pm w_{2}), \dots , (\pm x_{n}, \pm w_{n})\},$$
constructed from an antipodal free subset $\{(x_{1}, w_{1}),(x_{2}, w_{2}), \dots , (x_{n}, w_{n})\}$ of $S^{m}\times S^{M}$ contains exactly $4n$ distinct points.\ Applying Proposition \ref{tec1} to this set of points leads to the following outcome: the system
 $$
\sum_{\mu =1}^{n} \sum_{(i,j) \in I} c_{\mu}^{i,j}P_{k}^{m}((-1)^i x_{\mu}\cdot x)P_{l}^{M}((-1)^j w_{\mu} \cdot w)=0, \quad x \in S^M,\quad w \in  S^M, \quad (k,l) \in J_K,
$$
has just one solution, the zero one.\ Condition $(ii)$ is nothing but a re-writing of this fact.\ Conversely,  let $(x_{1}, w_{1}),(x_{2}, w_{2}), \dots , (x_{n}, w_{n})$ be distinct points of $S^{m}\times S^{M}$.\ We can employ Lemma \ref{antipodal} to pick an antipodal free subset
$\{(x'_{1}, w'_{1}),(x'_{2}, w'_{2}), \dots , (x'_{p}, w'_{p})\}$ of $S^{m}\times S^{M}$ so that
$$
\{(x_{\mu}, w_{\mu} ) :\mu=1,2,\ldots,n\} \subset \{ (\pm x'_{1}, \pm w'_{1}),(\pm x'_{2}, \pm w'_{2}), \dots , (\pm x'_{p}, \pm w'_{p})  \}. $$
If the system
$$
\sum_{\mu=1}^{n}c_{\mu}P_{k}^{m}(x_{\mu}\cdot x)P_{l}^{M}(w_{\mu}\cdot w)=0, \quad x \in S^m, \quad w \in S^M, \quad (k,l) \in J_K$$
has a nontrivial solution $(c_1, \ldots,c_n)$, it is promptly seen that the same will be true for the system
$$
\sum_{\mu =1}^{p}\sum_{(i,j) \in I} c_{\mu}^{i,j}P_{k}^{m}((-1)^i x'_{\mu}\cdot x)P_{l}^{M}((-1)^j w'_{\mu} \cdot w)=0, \quad x \in S^M, \quad w \in  S^M, \quad (k,l) \in J_K,
$$
that is, for the system described in $(ii)$.\eop\vspace*{3mm}

We now bring into play the disjoint union decomposition for $J_K$ explained in Section 1.\ In doing so, the previous theorem can be specialized to the following one.

\begin{thm} \label{maincarac} Let $K$ be a real, continuous, isotropic and positive definite kernel on $S^m\times S^M$.\ The following assertions are equivalent:\\
$(i)$ $K$ is strictly positive definite on $S^{m}\times S^{M}$;\\
$(ii)$ If $n$ is a positive integer and $\{(x_{1}, w_{1}),(x_{2}, w_{2}), \dots , (x_{n}, w_{n})\}$ is an antipodal free subset of $S^{m}\times S^{M}$, then the only solution $(d_{1}^{0,0},\ldots,d_{n}^{0,0},d_{1}^{1,0},\ldots,d_{n}^{1,0},
d_{1}^{0,1}, \ldots, d_{n}^{0,1},d_{1}^{1,1}, \ldots, d_{n}^{1,1})$ to the (block) system of equations
$$ \sum_{\mu=1}^{n} d_{\mu}^{i,j}P_{k}^{m}(x_{\mu}\cdot x)P_{l}^{M}(w_{\mu} \cdot w)=0, \quad x \in S^M, \quad w \in S^{M},\quad (k,l) \in J_K^{i,j}; \quad \quad (i,j) \in I,$$
is the zero vector.
\end{thm}
\pf Assume $K$ is strictly positive definite on $S^{m}\times S^{M}$.\ If $(ii)$ were not true, it would be possible to find an antipodal free subset $\{(x_{1}, w_{1}),(x_{2}, w_{2}), \dots,$ $(x_{n}, w_{n})\}$ of $S^m \times S^M$ for which the system described in $(ii)$ has a nontrivial solution.\ But then the subsystem
$$ \sum_{\mu=1}^{n} d_{\mu}^{i,j}P_{k}^{m}(x_{\mu}\cdot x)P_{l}^{M}(w_{\mu} \cdot w)=0, \quad (x,w) \in S^{m-1}\times S^{M-1}, \quad (k,l) \in J_K^{i,j},$$
would have a nontrivial solution $(d_1^{i,j},d_2^{i,j}, \ldots, d_n^{i,j})$, for at least one pair $(i,j) \in I$, say $(i,j)=(0,0)$ (in the other cases, the procedure is analogous).\ Now observe that for every $\mu \in \{1,2,\ldots, n\}$, the system of 4 linear equations
$$
\left\{\begin{array}{c}
  c_{\mu}^{0,0} + c_{\mu}^{1,0} + c_{\mu}^{0,1} + c_{\mu}^{1,1}= d_{\mu}^{0,0}\\
 c_{\mu}^{0,0}  -c_{\mu}^{1,0} + c_{\mu}^{0,1}  -c_{\mu}^{1,1}=0\\
 c_{\mu}^{0,0}  +c_{\mu}^{1,0}  -c_{\mu}^{0,1}  -c_{\mu}^{1,1}=0\\
 c_{\mu}^{0,0}  -c_{\mu}^{1,0}  -c_{\mu}^{0,1} + c_{\mu}^{1,1}=0
\end{array}\right.
$$
has a unique solution $(c_\mu^{0,0},c_\mu^{1,0},c_\mu^{0,1}, c_\mu^{1,1})$.\ If at least one $d_\mu^{0,0}$ is nonzero, then at least one of the solutions $(c_\mu^{0,0},c_\mu^{1,0},c_\mu^{0,1}, c_\mu^{1,1})$, $\mu \in \{1,2,\ldots,n\}$, is nontrivial.\
Therefore, returning to our previous arguments, the corresponding system in Proposition \ref{caracJean1}-$(ii)$ would have a nontrivial solution $$(c_{1}^{0,0},\ldots,c_{n}^{0,0},c_{1}^{1,0},\ldots,c_{n}^{1,0},
 c_{1}^{0,1}, \ldots, c_{n}^{0,1},c_{1}^{1,1}, \ldots, c_{n}^{1,1}).$$
Therefore, $K$ would not be strictly positive definite.\ Thus, $(i)$ implies $(ii)$.\ Conversely, let $\{(x_{1}, w_{1}),(x_{2}, w_{2}),\dots , (x_{n}, w_{n})\}$ be an antipodal free subset of $S^{m}\times S^{M}$ and consider the system described in $(ii)$.\ The sum $\sum_{(i,j) \in I}(-1)^{ik+jl}c_{\mu}^{i,j}$ remains constant in each $J_K^{i,j}$.\ As a matter of fact, it assumes the following values:
$$ \begin{array}{c}
d_\mu^{0,0}:=c_{\mu}^{0,0} + c_{\mu}^{1,0} + c_{\mu}^{0,1} + c_{\mu}^{1,1}, \quad (k,l)\in J_{f}^{0,0},\quad  \mu\in\{1,2,\ldots,n\},\\
d_\mu^{1,0}:=c_{\mu}^{0,0}  -c_{\mu}^{1,0} + c_{\mu}^{0,1}  -c_{\mu}^{1,1} \quad (k,l)\in J_{f}^{1,0}, \quad \mu\in\{1,2,\ldots,n\},\\
d_\mu^{0,1}:=c_{\mu}^{0,0}  +c_{\mu}^{1,0}  -c_{\mu}^{0,1}  -c_{\mu}^{1,1}, \quad (k,l) \in J_{f}^{0,1}, \quad \mu\in\{1,2,\ldots,n\},\\
d_\mu^{1,1}:=c_{\mu}^{0,0}  -c_{\mu}^{1,0}  -c_{\mu}^{0,1} + c_{\mu}^{1,1}, \quad (k,l) \in J_{f}^{1,1}. \quad \mu\in\{1,2,\ldots,n\}.
\end{array}
$$
If $(ii)$ holds, we can conclude that $d_\mu^{0,0}=d_\mu^{1,0}=d_\mu^{0,1}=d_\mu^{1,1}=0$.\ But that corresponds to
$$(c_{1}^{0,0},\ldots,c_{n}^{0,0},c_{1}^{1,0},\ldots,c_{n}^{1,0},
 c_{1}^{0,1}, \ldots, c_{n}^{0,1},c_{1}^{1,1}, \ldots, c_{n}^{1,1})=0.$$
In other words, the system in Proposition \ref{maincarac}-$(ii)$ has just one solution, the zero one.\ Thus $(i)$ holds. \eop

Before proving Theorem \ref{main2}, we need to recast a last technical result on the characterization of the strict positive definiteness of a real, continuous, isotropic and positive definite kernels on a single sphere (\cite{chen}).

\begin{lem} \label{CaracScho} Let $K$ be a real, continuous, isotropic and positive definite kernel on $S^m$ and consider the Fourier-Gegenbauer series representation for the isotropic part $K_r$ of $K$:
$$K_r=\sum_{k=0}^\infty a_k P_k^m ,\quad a_k\geq 0,\quad \sum_{k=0}^\infty a_k P_k^m(1)<\infty.$$
For distinct points $x_1, x_2, \ldots, x_n$ on $S^m$ and real numbers $c_1, c_2, \ldots, c_n$, the following assertions are equivalent:\\
$(i)$ $\sum_{\mu, \nu =1}^{n} c_{\mu}c_{\nu} g(x_{\mu} \cdot x_{\nu}) =0$;\\
$(ii)$ It holds
$$
\sum_{\mu=1}^{n}c_{\mu}P_{k}^{m}(x_{\mu}\cdot x)=0, \quad k \in \{k: a_k>0\}, \quad x \in S^{m-1}.$$
\end{lem}

\begin{thm} Let $K$ be a real, continuous, isotropic and positive definite kernel on $S^m\times S^M$.\ If it is strictly positive definite, then for each pair $(i,j) \in I$, there exists a sequence $\{(k_r^{i,j},l_r^{i,j})\}$ in $J_K^{i,j}$ so that $\lim_{r \to \infty}k_r^{i,j}=\lim_{r\to \infty}l_r^{i,j}=\infty$.
\end{thm}

\pf Assume there exists a positive integer $k_0$ so that $\{k: (k,l) \in J_K^{0,0}\} \subset \{0,1,\ldots,k_0\}$.\ For a positive odd integer $n$, let us define $(x_\mu,w_\mu) \in S^m \times S^M$
through the following expressions:
$$x_{\mu} =( \cos (2\pi \mu/n) , \sin (2\pi \mu/n),0, \dots ,0), \quad \mu=1,2,\ldots,n,$$
$$w_{\mu} =( \cos (2\pi \mu/n) , \sin (2\pi \mu/n),0, \dots ,0), \quad \mu=1,2,\ldots,n.$$
It is an obvious matter to certify that the set (of $n^2$ points)
$$\Gamma_n:=\{(x_\mu,w_\nu): \mu,\nu=1,2,\ldots,n\}$$
is antipodal free.\ Next, we define
$$c_{\mu} = (-1)^{\mu} ( e^{i\pi\mu / n} + e^{- i\pi\mu/ n} ), \quad \mu=1,2,\ldots,n,$$ and calculate the quadratic form
$$
QF(m):=\sum_{\mu, \nu=1}^{n} c_{\mu}c_{\nu}\sum_{k=0}^{k_0} P^{m}_{2k}(x_{\mu} \cdot x_{\nu}).
$$
Introducing Tchebyshev polynomials, we have that
\begin{eqnarray*}
QF(m) & = & \sum_{k=0}^{k_0}\sum_{j=0}^{k}c_{2k}^m(j)\sum_{\mu,\nu=1}^{n} c_{\mu}c_{\nu} \cos [2\pi (2k -2j)(\mu - \nu)/n ]\\
& = & \sum_{k=0}^{k_0}\sum_{j=0}^{k}c_{2k}^m(j)\left|\sum_{\mu=1}^{n} (-1)^{\mu} [ e^{i\pi\mu / n} + e^{- i\pi\mu/ n} ] e^{i2\pi \mu( 2k-2j)/n} \right|^2\\
& = &\sum_{k=0}^{k_0}\sum_{j=0}^{k}c_{2k}^m(j)\left|\sum_{\mu=1}^{n}  e^{i2\pi \mu ( 2k-2j-(n+1)/2)/n}  + \sum_{\mu=1}^{n}  e^{i2\pi \mu (  2k-2j-(n-1)/2)/n}\right|^2
\end{eqnarray*}
If $n> 4k_0+1$, then all the numbers $2k-2j-(n+1)/2$ and $2k-2j-(n-1)/2$ do not belong to $n\mathbb{Z}$.\ As so, the exponentials appearing above vanish and $QF(m)=0$.
An appeal to Lemma \ref{CaracScho} yields
$$\sum_{\mu=1}^{n}c_{\mu}P_{2k}^{m}(x_{\mu}\cdot x)=0, \quad k=0,1,\ldots,k_0,\quad  x \in S^m,$$
as long as $n> 4k_0+1$.\ Defining $d_{\mu,\nu}^{0,0}=c_\mu c_\nu$, $\mu,\nu=1,2,\ldots,n$, we finally deduce that
$$
\sum_{\mu , \nu=1}^{n} d_{\mu, \nu}^{0,0}P_{k}^{m}(x_{\mu}\cdot x)P_{l}^{M}(w_{\nu} \cdot w)=\sum_{\mu=1}^{n}c_{\mu}P_{2k}^{m}(x_{\mu}\cdot x)\sum_{\nu=1}^{n}c_{\nu}P_{2k}^{M}(w_{\nu}\cdot w)=0,
$$
whenever $(x,w) \in S^m \times S^M$, $(k,l) \in J_K^{0,0}$ and $n> 4k_0+1$.\ This information reveals that assertion $(ii)$ in Theorem \ref{maincarac}
does not hold for the set $\Gamma_n$ and $(i,j)=(0,0)$.\ Thus, $K$ cannot be strictly positive definite on $S^m \times S^M$.\ A similar procedure leads to the same conclusion, if we assume the existence of a positive integer $l_0$ such that $\{l: (k,l) \in J_K^{0,0}\} \subset \{0,1,\ldots, l_0\}$.\ Finally, the very same procedure can be adapted to hold in the case we replace the set $J_K^{0,0}$ with $J_K^{i,j}$, $(i,j) \in I\setminus\{(0,0)\}$.\ The proof is complete.\eop

\begin{thm} Let $K$ be a real, continuous, isotropic and positive definite kernel on $S^m\times S^M$.\ If for each pair $(i,j) \in I$, there exists a sequence $\{(k_r^{i,j},l_r^{i,j})\}$ in $J_K^{i,j}$ so that $\lim_{r \to \infty}k_r^{i,j}=\lim_{r\to \infty}l_r^{i,j}=\infty$, then $K$ is strictly positive definite on $S^m \times S^M$.
\end{thm}
\pf We intend to use Theorem \ref{maincarac}.\ So, for each pair $(i,j) \in I$, assume there exists a sequence $\{(k_r^{i,j},l_r^{i,j})\}$ in $J_K^{i,j}$ so that $\lim_{r \to \infty}k_r^{i,j}=\lim_{r\to \infty}l_r^{i,j}=\infty$.\  Let $\{(x_{1}, w_{1}),(x_{2}, w_{2}),$ $ \dots , (x_{n}, w_{n})\}$ be an antipodal free subset of $S^{m}\times S^{M}$ and suppose that
$$ \sum_{\mu=1}^{n} d_{\mu}^{i,j}P_{k}^{m}(x_{\mu}\cdot x)P_{l}^{M}(w_{\mu} \cdot w)=0, \quad (x, w) \in S^{m}\times S^{M},\quad (k,l) \in J_K^{i,j}; \quad (i,j) \in I.$$
For $\alpha \in \{1,2,\ldots,n\}$ fixed, let us choose $x=x_\alpha$ and $w=w_\alpha$ in the above system.\ For each pair $(i,j)$, we can plug in the appropriated sequence guaranteed by our assumption in each one of the four blocks of the system to obtain
$$d_{\alpha}^{i,j}P_{k_{r}^{i,j}}^{m}(1)P_{l_{r}^{i,j}}^{M}(1)+  \sum_{\mu \neq \alpha } d_{\mu }^{i,j}P_{k_{r}^{i,j}}^{m}(x_{\mu}\cdot x_{\alpha})P_{l_{r}^{i,j}}^{M}(w_{\mu} \cdot w_{\alpha})=0, \quad (i,j) \in I,$$
that is,
$$d_{\alpha}^{i,j}+  \sum_{\mu \neq \alpha } d_{\mu }^{i,j}R_{k_{r}^{i,j}}^{m}(x_{\mu}\cdot x_{\alpha})R_{l_{r}^{i,j}}^{M}(w_{\mu} \cdot w_{\alpha})=0,\quad (i,j) \in I.$$
Taking into account that $x_{\mu}\cdot x_{\alpha} \neq \pm 1$, $\mu \neq \alpha$, that $w_{\mu} \cdot w_{\alpha} \neq \pm1$, $\mu \neq \alpha$, and Lemma \ref{gege}-$(ii)$, we can let $r \to \infty$ in each of the four resulting equations to deduce that $d_{\alpha}^{i,j}=0$.\eop

Recalling the remark at the end of Section 3, it follows that Theorem \ref{main2} holds when either $m=\infty$ or $M=\infty$.\ A version of the theorem in the case in which one of the spheres is a circle but the other one is not, is still an open question.


%
%

\vspace*{5mm}

\noindent
Departamento de
Matem\'atica,\\ ICMC-USP - S\~ao Carlos, Caixa Postal 668,\\
13560-970 S\~ao Carlos SP, Brasil\\ e-mails: jeanguella@gmail.com; menegatt@icmc.usp.br


\begin{thebibliography}{1}

\bibitem{askey} Askey, R., Orthogonal polynomials and special functions.\ Society for Industrial and Applied Mathematics, Philadelphia, Pa., 1975.
\bibitem{atkinson} Atkinson, K.; Han, Weimin, Spherical harmonics and approximations on the unit sphere: an introduction.\ Lecture Notes in Mathematics, 2044. Springer, Heidelberg, 2012.
\bibitem{berg} Berg, C.; Christensen, J. P. R.; Ressel, P., Harmonic analysis on semigroups.\ Theory of positive definite and related functions.\ Graduate Texts in Mathematics, 100.\ Springer-Verlag, New York, 1984.
\bibitem{bingham} Bingham, N. H., Positive definite functions on spheres.\ {\em Proc. Cambridge Philos. Soc.} 73 (1973), 145-156.
\bibitem{chen} Chen, Debao; Menegatto, V. A.; Sun, Xingping, A necessary and sufficient condition for strictly positive definite functions on spheres.\ {\em Proc. Amer. Math. Soc.}
 131 (2003), no. 9, 2733-2740.
 \bibitem{cheney} Cheney, E. W., Approximation using positive definite functions.\ Approximation theory VIII, Vol. 1 (College Station, TX, 1995), 145-168, Ser. Approx. Decompos., 6, World Sci. Publ., River Edge, NJ, 1995.
\bibitem{dai} Dai, Feng; Xu, Yuan, Approximation theory and harmonic analysis on spheres and balls.\ Springer Monographs in Mathematics.\ Springer, New York, 2013.
\bibitem{gneiting} Gneiting, T., Strictly and non-strictly positive definite functions on spheres.\ {\em Bernoulli} 19 (2013), no. 4, 1327-1349.
\bibitem{groemer} Groemer, H., Geometric applications of Fourier series and spherical harmonics.\ Encyclopedia of Mathematics and its Applications, 61.\ Cambridge University Press, Cambridge, 1996.
\bibitem{jean} Guella, J. C.; Menegatto, V. A.; Peron, A. P., An extension of a theorem of Schoenberg to products of spheres.\ ArXiv:1503.08174v1.
\bibitem{guella}  Guella, J. C.; Menegatto, V. A.; Peron, A. P., Strictly positive definite kernels on a product of circles.\ Arxiv:1505.01169.
\bibitem{muller} M{\"u}ller, C., Analysis of spherical symmetries in Euclidean spaces.\ Applied Mathematical Sciences, 129.\ Springer-Verlag, New York, 1998.
\bibitem{schoen} Schoenberg, I. J.; Positive definite functions on spheres.\ {\em Duke Math. J.} 9, (1942), 96-108.
\bibitem{szego}
Szeg\"{o}, G., Orthogonal polynomials.\ Fourth edition.\ American Mathematical Society, Colloquium Publications, Vol. XXIII.\ American Mathematical Society, Providence, R.I., 1975.





\end{thebibliography}
\end{document}